\documentclass[12pt]{article}
\usepackage{latexsym, amssymb, amsmath, amscd, amsfonts, epsfig, graphicx, colordvi,verbatim,ifpdf}
\usepackage{amsfonts, amsmath, amssymb,extarrows}
\usepackage{amssymb,amsfonts,amsmath,latexsym,epsfig,cite, psfrag,eepic,color}
\usepackage{amscd,graphics}
\usepackage{latexsym, amssymb,  amsmath,amscd, amsfonts, epsfig, graphicx, colordvi,amsthm}

\usepackage{graphicx}
\usepackage{color}
\usepackage{ifpdf}
\usepackage{fancybox}

\usepackage{float}

\newtheorem{thm}{Theorem}[section]

\newtheorem{conj}[thm]{Conjecture}

\newtheorem{defi}[thm]{Definition}
\newtheorem{lem}[thm]{Lemma}
\newtheorem{core}[thm]{Corollary}

\def\pf{\noindent{\it Proof.} }
\def\qed{\nopagebreak\hfill{\rule{4pt}{7pt}}
\medbreak}

\numberwithin{equation}{section}

\def\qed{\nopagebreak\hfill{\rule{4pt}{7pt}}
\medbreak}

\newlength{\boxedparwidth}
  {\begin{center} \begin{tabular}{|@{\hspace{.315in}}c@{\hspace{.15in}}|}
                  \hline \\ \begin{minipage}[t]{\boxedparwidth}
                  \setlength{\parindent}{.25in}}%
  {\end{minipage} \\ \\ \hline \end{tabular} \end{center}}

\allowdisplaybreaks
\begin{document}

\begin{center}

 {\large \bf Monotonicity properties for ranks of overpartitions}
\end{center}
 
\vskip 5mm

\begin{center}
{  Huan Xiong}$^{1}$ 
and {Wenston J.T. Zang}$^{2}$ \vskip 2mm

 $^{1}$Institut de Recherche Math\'ematique Avanc\'ee, UMR 7501 \\ 
 Universit\'e de Strasbourg et CNRS,  F-67000 Strasbourg, France\\[6pt]
    $^{2}$Institute of Advanced Study of Mathematics \\
   Harbin Institute of Technology, Heilongjiang 150001, P.R. China\\[6pt]

   \vskip 2mm

Email:     $^1$xiong@math.unistra.fr, $^2$zang@hit.edu.cn 
\end{center}

\vskip 6mm \noindent {\bf Abstract.} The rank of partitions play an important role in the combinatorial interpretations of several Ramanujan's famous congruence formulas. In 2005 and 2008, the $D$-rank and $M_2$-rank of an overpartition were introduced by Lovejoy,  respectively. Let $\overline{N}(m,n)$ and $\overline{N2}(m,n)$ denote the number of overpartitions of $n$ with $D$-rank $m$ and $M_2$-rank $m$, respectively. In 2014, Chan and Mao proposed a conjecture on monotonicity properties of $\overline{N}(m,n)$ and $\overline{N2}(m,n)$. In this paper, we prove the Chan-Mao monotonicity conjecture. To be specific, we show that for any integer $m$ and nonnegative integer $n$,  $\overline{N2}(m,n)\leq \overline{N2}(m,n+1)$; and for $(m,n)\neq (0,4)$ with $n\neq\, |m| +2$, we have $\overline{N}(m,n)\leq \overline{N}(m,n+1)$. Furthermore, when $m$ increases, we prove that $\overline{N}(m,n)\geq \overline{N}(m+2,n)$ and $\overline{N2}(m,n)\geq \overline{N2}(m+2,n)$ for any $m,n\geq 0$, which is an analogue of Chan and Mao's result for partitions.



\section{Introduction}


The aim of this paper is to study monotonicity properties of the $D$-rank and $M_2$-rank on overpartitions and therefore prove a conjecture of Chan and Mao \cite{Chan-Mao-2014}.

Recall that a partition of a nonnegative integer $n$ is a finite weakly decreasing sequence of positive integers $\lambda = (\lambda_1, \lambda_2, \ldots, \lambda_\ell)$ with $\sum_{1\leq i\leq \ell}\lambda_i=n$. Here $\lambda_1, \lambda_2, \ldots, \lambda_\ell$ are called parts of the partition $\lambda$ (see \cite{Andrews-1976}). The rank of a partition was defined by Dyson \cite{Dyson-1944} as the largest part of the partition minus the number of parts.  Dyson first conjectured and then proved by Atkin and Swinnerton-Dyer \cite{Atkin-Swinnerton-Dyer-1954} that the rank can provide combinatorial interpretations to the following Ramanujan's famous congruence for the partition function modulo $5$ and $7$, respectively: 
   \begin{eqnarray}
 p(5n+4)&\equiv& 0 \pmod 5, \label{con-5}\\[3pt]
p(7n+5) &\equiv& 0 \pmod 7, \label{con-7}
\end{eqnarray}
where $p(n)$ denotes the number of partitions of $n$. Since then, various results on the rank of partitions have been obtained by many mathematicians (For example, see \cite{Andrews-Garvan-1988,Andrews-Chan-Kim-2013,Andrews-Lewis-2000,Atkin-Garvan-2003,Bringmann-2009,Bringmann-Kane-2010,Bringmann-Mahlburg-2009,Bringmann-Mahlburg-2014,Bringmann-Mahlburg-2011,Bringmann-Lovejoy-2007,Chan-Mao-2014,Chen-Ji-Zang-2015,Chen-Ji-Zang-2016,Garvan-1990,Garvan-Stanton-1990,Kane-2004,Lewis-1991-1,Lewis-1991-2,Lewis-1992-1,Lewis-1992-2,Lewis-1997,Lewis-2009,Lewis-1994,Santa-1992}).

Let $N(m,n)$ denote the number of partitions of $n$ with rank $m$. Chan and Mao \cite{Chan-Mao-2014} established the following monotonicity properties for $N(m,n)$.

\begin{thm}[Chan and Mao \cite{Chan-Mao-2014}]\label{thm-cm-1}
For $n\geq 12$, $m\geq 0$ and $n\neq m+2$,
\begin{equation}\label{thm-cm-1-eq}
N(m,n)\geq N(m,n-1).
\end{equation}
\end{thm}

\begin{thm}[Chan and Mao \cite{Chan-Mao-2014}]\label{thm-cm-2}
For $n\geq 0$ and $m\geq 0$,
\begin{equation}\label{thm-cm-2-eq}
N(m,n)\geq N(m+2,n).
\end{equation}
\end{thm}

At the end of their paper, Chan and Mao \cite{Chan-Mao-2014} proposed a conjecture on monotonicity properties of the $D$-rank and $M_2$-rank of an overpartition. Recall that an overpartition was defined by Corteel and Lovejoy \cite{Corteel-Lovejoy-2004}  as a partition of $n$ in
which the first occurrence of a part may be overlined. For example, there are $14$ overpartitions of $4$:
\[\begin{array}{lllllllllllll}
(4),&(\bar{4}),&(3,1),&(\bar{3},1),&(3,\bar{1}),&(\bar{3},\bar{1}),&(2,2),\\[5pt]
(\bar{2},2)&(2,1,1),&(\bar{2},1,1),&(2,\bar{1},1),&(\bar{2},\bar{1},1),& (1,1,1,1),&(\bar{1},1,1,1).
\end{array}
\]

Lovejoy \cite{Lovejoy-2005} defined the $D$-rank of an overpartition as the largest part minus the number of parts, which is an analogue of the rank on ordinary partitions. Let $\overline{N}(m,n)$ denote the number of overpartitions of $n$ with $D$-rank $m$. Lovejoy \cite[Proposition 1.1]{Lovejoy-2005} gave the following generating function of $\overline{N}(m,n)$:
\begin{equation}\label{equ-gf-over-rank}
\sum_{n= 0}^\infty\sum_{m=-\infty}^\infty \overline{N}(m,n) z^mq^n=\sum_{k=0}^\infty\frac{(-1;q)_k \, q^{k(k+1)/2}}{(zq;q)_k(q/z;q)_k}.
\end{equation}
Here and throughout the rest of this paper, we adopt the common $q$-series notation \cite{Andrews-1976}:
\begin{align*}
(a;q)_\infty=\prod_{n=0}^\infty(1-aq^n) \quad \text{and} \quad
(a;q)_n&=\frac{(a;q)_\infty}{(aq^n;q)_\infty}.
\end{align*}

 The $M_2$-rank on overpartitions was also introduced by Lovejoy \cite{Lovejoy-2008}. For an overpartition $\lambda$, let $\lambda_1$ denote the largest part of $\lambda$, $\ell(\lambda)$ denote the number of parts of $\lambda$, and $\lambda_o$ denote the partition consisting of the non-overlined odd parts of $\lambda$. Then define
\begin{equation}
M_2\text{-rank}(\lambda)=\left\lfloor\frac{\lambda_1}{2}\right\rfloor-\ell(\lambda)
+\ell(\lambda_o)-\chi(\lambda),
\end{equation}
where $\chi(\lambda)=1$ if the largest part of $\lambda$ is odd and non-overlined, and otherwise $\chi(\lambda)=0$.

For instance, let $\lambda=(\overline{7},5,\overline{4},4,\overline{2},2,1,1)$. Then $\lambda_1=7$, $\ell(\lambda)=8$, $\lambda_o=(5,1,1)$, $\ell(\lambda_o)=3$ and $\chi(\lambda)=0$. Therefore,
$$
M_2\text{-rank}(\lambda)=3-8
+3=-2.
$$

Let $\overline{N2}(m,n)$ denote the number of overpartitions of $n$ with $M_2$-rank $m$.  Lovejoy \cite{Lovejoy-2008} found the generating function of $\overline{N2}(m,n)$ as follows:
\begin{equation}\label{equ-gf-over-m2-rank}
\sum_{n=0}^\infty\sum_{m=-\infty}^\infty \overline{N2}(m,n) z^mq^n=
\sum_{k=0}^\infty\frac{(-1;q)_{2k}q^{k}}{(zq^2;q^2)_k(q^2/z;q^2)_k}.
\end{equation}

Various results on the $D$-rank and $M_2$-rank of overpartitions can be found in \cite{Andrews-Chan-Kim-Osburn-2016,Andrews-Dixit-Schultz-Yee-2016,Bringmann-Lovejoy-Osburn-2009,Garvan-Jennings-2014,Jennings-2015-1,Jennings-2015-2,Jennings-2016,Lovejoy-2005,Lovejoy-2008,Lovejoy-Osburn-2008,Lovejoy-Osburn-2010}.
In 2014, Chan and Mao \cite{Chan-Mao-2014} proposed the following monotonicity conjecture on $\overline{N}(m,n)$ and $\overline{N2}(m,n)$:

\begin{conj}[Chan and Mao \cite{Chan-Mao-2014}]\label{main-conj}
For $(m,n)\neq (0,4)$ with $n\neq |m|+2$, we have
\begin{equation}\label{ine-nmn-ge-nmn-1}
\overline{N}(m,n)\geq \overline{N}(m,n-1).
\end{equation}
For $m \in \mathbb{Z}$ and $n\geq 0$,
\begin{equation}\label{ine-n2mn-ge-n2mn-1}
\overline{N2}(m,n)\geq \overline{N2}(m,n-1).
\end{equation}
\end{conj}

The main purpose of this paper is to give analogues of Theorems  \ref{thm-cm-1} and \ref{thm-cm-2}. To be specific, we obtain the following results:

\begin{thm}\label{main-thm-1}
For $m,n\geq 0$ with $n\neq m+2$ and $(m,n)\neq (0,4)$,
\begin{equation}\label{ine-nmn-ge-nmn-1-1-1}
\overline{N}(m,n)\geq \overline{N}(m,n-1).
\end{equation}
For $m,n\geq 0$, we have
\begin{equation}\label{ine-nmn-ge-nmn-1-1-2}
\overline{N2}(m,n)\geq \overline{N2}(m,n-1).
\end{equation}
\end{thm}

\begin{thm}\label{main-thm-2}
For $m,n\geq 0$, we have
\begin{equation}\label{ine-nmn-ge-nmn-1-1-2-1}
\overline{N}(m,n)\geq \overline{N}(m+2,n),
\end{equation}
and
\begin{equation}\label{ine-nmn-ge-nmn-1-1-2-2}
\overline{N2}(m,n)\geq \overline{N2}(m+2,n),
\end{equation}
\end{thm}
By the generating functions \eqref{equ-gf-over-rank} and \eqref{equ-gf-over-m2-rank}, it is easy to see that $\overline{N}(-m,n)=\overline{N}(m,n)$ and $\overline{N2}(-m,n)=\overline{N2}(m,n)$. Therefore Theorem \ref{main-thm-1} verifies Conjecture \ref{main-conj}.

This paper is organized as follows. Some preliminary results are given in Section \ref{sec:2}. Then in Section \ref{sec:3}, we establish a nonnegativity result Lemma \ref{lem-1} and use it to give a proof of Theorem \ref{main-thm-1}. Section \ref{sec:4} is devoted to prove Theorem \ref{main-thm-2}.

\section{Preliminary}\label{sec:2}

In order to prove Theorems \ref{main-thm-1} and Theorem \ref{main-thm-2}, we need to recall the definition of a function $f_{m,k}(q)$, which was first given by Chan and Mao \cite{Chan-Mao-2014}.

\begin{defi}
Define $f_{m,k}(q)$ as coefficients in the following formal power series:
\begin{equation}\label{equ-gf-f-mk}
\sum_{m=-\infty}^\infty z^m f_{m,k}(q):=\frac{1-q}{(zq;q)_k(q/z;q)_k}.
\end{equation}
\end{defi}

When $k=0$, by definition we see that $f_{0,0}(q)=1-q$ and $f_{m,0}(q)=0$ for all $m\neq 0$. Chan and Mao \cite[Lemma 9]{Chan-Mao-2014} gave the following expressions for  $f_{m,1}(q)$ and $f_{m,2}(q)$.

\begin{thm}[Chan and Mao \cite{Chan-Mao-2014}]\label{lem-chan-mao-two-lem-1}
For all integer $m$,
\begin{equation}\label{equ-f-1-k-gf}
f_{m,1}(q)=\sum_{n=\mid m\mid}^\infty  (-1)^{m+n} q^n=\frac{q^{\mid m\mid}}{1+q}.
\end{equation}
For $m=0$,
\begin{eqnarray}\label{equ-f-2-k-gf}
f_{0,2}(q)
=-q+\frac{1}{1-q^3}+\frac{q^2}{1-q^4}+\frac{q^8}{(1-q^3)(1-q^4)},
\end{eqnarray}
and for $m\neq 0$,
\begin{eqnarray}\label{equ-f-2-k-gf-m-neq-0}
f_{m,2}(q)= q^{\mid m\mid}\left(\frac{1-q^{\mid m\mid+1}}{(1-q^2)(1-q^3)}+\frac{q^{\mid m\mid+3}}{(1-q^3)(1-q^4)}\right).
\end{eqnarray}
\end{thm}

 Chan and Mao \cite[Lemma 11]{Chan-Mao-2014} also found the following nonnegative property for $f_{m,k}(q)$ when $k\geq 2$. For the remainder part of this paper, let $\{b_n\}_{n=0}^\infty$ be any sequence of nonnegative integers but not necessarily the same in different equations.

\begin{thm}[Chan and Mao \cite{Chan-Mao-2014}]\label{lem-chan-mao-two-lem-2}
When $k\geq 2$,
\begin{eqnarray}
f_{0,k}(q)&=&-q+q^2+\sum_{n=0}^\infty b_nq^n;\label{equ-f0k-bn}\\[3pt]
f_{1,k}(q)&=& q^{k+2}+\sum_{n=0}^\infty b_nq^n;\label{equ-f1k-bn}\\[3pt]
f_{m,k}(q)&=& \sum_{n=0}^\infty b_nq^n,\quad\text{for }m\geq 2.\label{equ-fmk-bn}
\end{eqnarray}
\end{thm}

By definition, it is easy to check that the constant term of $f_{0,k}(q)$ is equal to $1$. Hence \eqref{equ-f0k-bn} yields the following corollary:

\begin{core}\label{coro-1}
When $k\geq 2$,
\begin{equation}
f_{0,k}(q)=1-q+q^2+\sum_{n=0}^\infty b_nq^n.
\end{equation}
\end{core}

We also need the following two lemmas in \cite{Chan-Mao-2014}. 

\begin{lem}[See Lemma 8 of  \cite{Chan-Mao-2014}]\label{lem-8}
When $k\geq 0$, we have
$$
f_{m,k+1}(q)=\sum_{n=-\infty}^\infty f_{n,k}(q)\,\frac{q^{{(k+1)}\mid m-n\mid}}{1-q^{2k+2}}.
$$
\end{lem}

\begin{lem}[See Lemma 10 of 
 \cite{Chan-Mao-2014}]\label{lem-10}
For any positive integer $m$,
$$
\frac{1-q^{m+1}}{(1-q^2)(1-q^3)}
$$
has nonnegative power series coefficients.
\end{lem}


\section{The proof of Theorem \ref{main-thm-1}} \label{sec:3}

In this section, we give a proof of Theorem \ref{main-thm-1}. To this end, we need the following lemma.

\begin{lem}\label{lem-1}
For any nonnegative integer $a,b$ and $c$, the coefficient of $q^n$ in
\[\frac{q^a}{1+q^c}+\frac{q^b}{(1-q^3)(1-q^4)}\]
is nonnegative for $n\geq b+6$.
\end{lem}

\pf It is clear that
\[\frac{q^b}{(1-q^3)(1-q^4)}=\sum_{i=0}^\infty\sum_{j=0}^\infty \, q^{b+3i+4j}.\]
Note that for any $n\geq 6$, there exists $i,j\geq 0$ such that $3i+4j=n$. To be specific,
\begin{equation}
(i,j)=
\begin{cases}
(k,0) & \text{if }n=3k;\\[3pt]
(k-1,1)& \text{if }n=3k+1;\\[3pt]
(k-2,2) & \text{if }n=3k+2.
\end{cases}
\end{equation}
 Hence we see that, the coefficient of $q^n$ in
\begin{equation}
\frac{q^b}{(1-q^3)(1-q^4)}
\end{equation}
is at least $1$.
On the other hand,
\begin{align}
\frac{q^a}{1+q^c}&=\sum_{m=0}^\infty (-1)^{m}q^{cm+a}.
\end{align}
Evidently, for any nonnegative integer $n$, the coefficient of $q^n$ in $\sum_{m=0}^\infty (-1)^{m}q^{cm+a}$ is either $-1$, $0$ or $1$. Thus when $n\geq b+6$, the coefficient of $q^n$ in
\[\frac{q^a}{1+q^c}+\frac{q^b}{(1-q^3)(1-q^4)}\]
is nonnegative. This yields the  desired result.\qed

We are now in a position to prove Theorem \ref{main-thm-1}.

{\noindent \it Proof of Theorem \ref{main-thm-1}.}
We first prove \eqref{ine-nmn-ge-nmn-1-1-1} with the aid of Lemma \ref{lem-1}, and then show \eqref{ine-nmn-ge-nmn-1-1-2}.

 From \eqref{equ-gf-over-rank}, it is clear to see that
\begin{equation}\label{equ-dif-over-rank}
1+\sum_{n=1}^\infty\sum_{m=-\infty}^\infty \left(\,\overline{N}(m,n)-\overline{N}(m,n-1)\,\right) z^mq^n=\sum_{k=0}^\infty\frac{(-1;q)_k \, q^{k(k+1)/2}(1-q)}{(zq;q)_k(q/z;q)_k}.
\end{equation}
By the definition of $f_{m,k}(q)$ (see \eqref{equ-gf-f-mk}), we derive that
\begin{equation}\label{equ-over-rank-nmn-fmk}
1+\sum_{n=1}^\infty\sum_{m=-\infty}^\infty \left(\overline{N}(m,n)-\overline{N}(m,n-1)\right)  z^mq^n=
\sum_{m=-\infty}^\infty z^m\sum_{k=0}^\infty (-1;q)_k \, q^{k(k+1)/2}f_{m,k}(q).
\end{equation}
Hence for fixed integer $m\neq 0$,
\begin{equation}\label{equ-over-rank-nmn-nmn1-1}
\sum_{n=1}^\infty \left(\,\overline{N}(m,n)-\overline{N}(m,n-1)\,\right)  q^n=\sum_{k=0}^\infty (-1;q)_k \, q^{k(k+1)/2}f_{m,k}(q).
\end{equation}
When $m=0$,
by \eqref{equ-over-rank-nmn-fmk}, 
  \eqref{equ-over-rank-nmn-nmn1-1} and Theorem \ref{lem-chan-mao-two-lem-1} we find that
\begin{align}\label{equ-over-rank-n0n-n0n1-0}
&\sum_{n=1}^\infty \left(\,\overline{N}(0,n)-\overline{N}(0,n-1)\,\right) q^n\nonumber\\[3pt]
=&\ -q+\frac{2q}{1+q}+2(1+q)q^3\left(-q+\frac{1}{1-q^3}+\frac{q^2}{1-q^4}
+\frac{q^8}{(1-q^3)(1-q^4)}\right)\nonumber\\[3pt]
&+\sum_{k=3}^\infty (-1;q)_k \, q^{k(k+1)/2}f_{0,k}(q).
\end{align}
By Corollary \ref{coro-1}, we derive that
\begin{align}\label{equ-over-rank-n0n-n0n1-1}
&\sum_{n=1}^\infty \left(\,\overline{N}(0,n)-\overline{N}(0,n-1)\,\right) q^n\nonumber\\[3pt]
=&-q-2q^4-2q^5+\frac{2(1+q)q^3}{1-q^3}+\frac{2(1+q)q^5}{1-q^4}
+\frac{2q^{12}}{(1-q^3)(1-q^4)}
\nonumber\\[3pt]
&+\frac{2q}{1+q}
+\frac{2q^{11}}{(1-q^3)(1-q^4)}\nonumber\\[3pt]
&+\sum_{k=3}^\infty (-1;q)_k \, q^{k(k+1)/2}\left(1-q+q^2+\sum_{n=0}^\infty b_nq^n\right).
\end{align}
The last term in \eqref{equ-over-rank-n0n-n0n1-1} can be transformed as follows:
\begin{align}\label{equ-over-rank-n0n-n0n1-1-term3}
&\sum_{k=3}^\infty (-1;q)_k \, q^{k(k+1)/2}\left(1-q+q^2+\sum_{n=0}^\infty b_nq^n\right)\nonumber\\[3pt]
=&
\sum_{k=3}^\infty 2(1+q)(-q^2;q)_{k-2}\, q^{k(k+1)/2}(1-q+q^2)+\sum_{k=3}^\infty(-1;q)_k \, q^{k(k+1)/2}\sum_{n=0}^\infty b_nq^n\nonumber\\[3pt]
=&\sum_{k=3}^\infty 2(1+q^3)(-q^2;q)_{k-2} q^{k(k+1)/2}+ \sum_{k=3}^\infty (-1;q)_k \, q^{k(k+1)/2}\sum_{n=0}^\infty b_nq^n,
\end{align}
which clearly has nonnegative coefficients. Moreover, by Lemma \ref{lem-1}, the coefficient of $q^n$ in
\[\frac{2q}{1+q}
+\frac{2q^{11}}{(1-q^3)(1-q^4)}\]
is nonnegative for $n\geq 17$. From the above analysis, we see that
\[\overline{N}(0,n)\geq \overline{N}(0,n-1)\]
for $n\geq 17$. It is trivial to check that for $1\leq n\leq 16$,
\[\overline{N}(0,n)\geq \overline{N}(0,n-1)\]
except for $n=2$ or $n=4$. Therefore Theorem \ref{main-thm-1} holds for $m=0$.

We now assume that $m\geq 1$.
Substituting \eqref{equ-f-1-k-gf} and \eqref{equ-f-2-k-gf-m-neq-0} into \eqref{equ-over-rank-nmn-nmn1-1}, we have
\begin{align}
&\sum_{n=1}^\infty \left(\,\overline{N}(m,n)-\overline{N}(m,n-1)\,\right)  q^n\nonumber\\[3pt]
=&\,\frac{2q^{m+1}}{1+q}+\sum_{k=3}^\infty (-1;q)_k \, q^{k(k+1)/2}f_{m,k}(q)\nonumber\\[3pt]
&+2(1+q)q^{m+3}\left(\frac{1-q^{m+1}}{(1-q^2)(1-q^3)}
+\frac{q^{m+3}}{(1-q^3)(1-q^4)}\right).
\end{align}
From Theorem \ref{lem-chan-mao-two-lem-2}, we see that for $k\geq 3$, $f_{m,k}(q)$ has nonnegative coefficients. We proceed to show  the coefficients of $q^n$ in
\begin{equation}\label{equ-rank-over-neq-1-3-1}
\frac{2q^{m+1}}{1+q}+2(1+q)q^{m+3}\left(\frac{1-q^{m+1}}{(1-q^2)(1-q^3)}+\frac{q^{m+3}}{(1-q^3)(1-q^4)}\right)
\end{equation}
is nonnegative for all $n\geq m+3$.

We first assume that $m\neq 1,3$. In this case, we transform \eqref{equ-rank-over-neq-1-3-1} as follows:
\begin{align}\label{equ-rank-over-neq-1-3-1-1-1}
&\frac{2\,q^{m+1}}{1+q}+2(1+q)q^{m+3}
\left(\frac{1-q^{m+1}}{(1-q^2)(1-q^3)}+\frac{q^{m+3}}{(1-q^3)(1-q^4)}\right)\nonumber\\[3pt]
=&\,\frac{2\,q^{m+1}}{1+q}+2q^{m+4}\frac{1-q^{m+1}}{(1-q^2)(1-q^3)}\nonumber\\[3pt]
&+2q^{m+3}\frac{1-q^{m+1}}{(1-q^2)(1-q^3)}+2(1+q)\frac{q^{2m+6}}{(1-q^3)(1-q^4)}.
\end{align}
By Lemma \ref{lem-10}, we find that
\[2q^{m+3}\frac{1-q^{m+1}}{(1-q^2)(1-q^3)}\]
has nonnegative coefficients in $q^n$ for all $n\geq 1$. Moreover,
\begin{eqnarray*}
\frac{2q^{m+1}}{1+q}+2q^{m+4}\frac{1-q^{m+1}}{(1-q^2)(1-q^3)}&=&\frac{2q^{m+1}}{1+q}+
2q^{m+4}\frac{1-q^3+q^3-q^{m+1}}{(1-q^2)(1-q^3)}\nonumber\\[3pt]
&=&\frac{2q^{m+1}}{1+q}+\frac{2q^{m+4}}{1-q^2}+2q^{m+7}\frac{1-q^{m-2}}{(1-q^2)(1-q^3)}
\nonumber\\[3pt]
&=&2q^{m+1}\frac{1-q+q^3}{1-q^2}+2q^{m+7}\frac{1-q^{m-2}}{(1-q^2)(1-q^3)}\nonumber\\[3pt]
&=&\frac{2q^{m+1}}{1-q^2}-2q^{m+2}+2q^{m+7}\frac{1-q^{m-2}}{(1-q^2)(1-q^3)}.
\end{eqnarray*}
Notice that when $m\neq 1,3$, by Lemma \ref{lem-10} we obtain
\[2q^{m+7}\frac{1-q^{m-2}}{(1-q^2)(1-q^3)}=\sum_{n=0}^\infty b_nq^n.\]
This yields that \eqref{equ-rank-over-neq-1-3-1-1-1} has nonnegative coefficients in $q^n$ for $n\geq m+3$, as desired.

It remains to consider the case $m=1$ or $3$.
For $m=1$, it is trivial to calculate that \eqref{equ-rank-over-neq-1-3-1} is equal to
\begin{equation}\label{equ-rank-over-neq-1-3-2}
\frac{2q^2}{1+q}+\frac{2q^4+2q^5}{(1-q^3)(1-q^4)}.
\end{equation}
From Lemma \ref{lem-1}, we see that for $n\geq 10$, the coefficient of $q^n$ in
\[\frac{2q^2}{1+q}+\frac{2q^4}{(1-q^3)(1-q^4)}\]
is nonnegative. Hence we derive that $\overline{N}(1,n)\geq N(1,n-1)$ for $n\geq 10$.
It is trivial to check that for $4\leq n\leq 9$, $\overline{N}(1,n)\geq N(1,n-1)$ also holds. This yields the case for $m=1$.

Finally, for $m=3$, \eqref{equ-rank-over-neq-1-3-1} is equal to:
\begin{equation}\label{equ-rank-over-neq-3-3-1}
\frac{2q^4}{1+q}+\frac{2q^{12}}{(1-q^3)(1-q^4)}+\frac{2q^{13}}{(1-q^3)(1-q^4)}+
\frac{2(1+q)(1+q^2)q^6}{1-q^3}.
\end{equation}
Using Lemma \ref{lem-1}, we find that for $n\geq 18$, the coefficient of $q^n$ in
\begin{equation}\label{equ-rank-over-neq-3-3-2}
\frac{2q^4}{1+q}+\frac{2q^{12}}{(1-q^3)(1-q^4)}
\end{equation}
is nonnegative. This yields that $\overline{N}(3,n)\geq \overline{N}(3,n-1)$ for $n\geq 18$. After checking $\overline{N}(3,n)\geq \overline{N}(3,n-1)$ for $6\leq n\leq 17$, we find that \eqref{ine-nmn-ge-nmn-1-1-1} is valid for $m=3$.

We next prove \eqref{ine-nmn-ge-nmn-1-1-2}. From \eqref{equ-gf-over-m2-rank}, we see that
\begin{align}
&\qquad 1+\sum_{n=1}^\infty\sum_{m=-\infty}^\infty \left(\,\overline{N2}(m,n)-\overline{N2}(m,n-1)\,\right) z^mq^n\nonumber\\[3pt]
&=(1-q)
\sum_{k=0}^\infty\frac{(-1;q)_{2k}q^{k}}{(zq^2;q^2)_k(q^2/z;q^2)_k}\nonumber\\[3pt]
&=1-q+2\sum_{k=1}^\infty\frac{(1-q^2)(-q^2;q)_{2k-2}q^k}{(zq^2;q^2)_k(q^2/z;q^2)_k}\nonumber\\[3pt]
&=1-q+2\sum_{k=1}^\infty(-q^2;q)_{2k-2}q^k\sum_{m=-\infty}^\infty z^m f_{m,k}(q^2).
\end{align}
Hence
\begin{equation}
1+\sum_{n=1}^\infty \left(\,\overline{N2}(0,n)-\overline{N2}(0,n-1)\,\right) q^n=1-q+2\sum_{k=1}^\infty(-q^2;q)_{2k-2}q^k f_{0,k}(q^2),
\end{equation}
and for $m\geq 1$,
\begin{equation}\label{equ-fg-n2-0-1}
\sum_{n=1}^\infty \left(\,\overline{N2}(m,n)-\overline{N2}(m,n-1)\,\right) q^n=2\sum_{k=1}^\infty(-q^2;q)_{2k-2}\,q^k f_{m,k}(q^2).
\end{equation}
Similar to the proof of \eqref{ine-nmn-ge-nmn-1-1-1}, we first assume that $m=0$. From Theorem \ref{lem-chan-mao-two-lem-1} and Corollary \ref{coro-1}, we deduce that
\begin{align}\label{equ-gf-n2mn-0-1}
&\qquad 1+\sum_{n=1}^\infty \left(\,\overline{N2}(0,n)-\overline{N2}(0,n-1)\,\right) q^n\nonumber\\[3pt]
&=1-q+
\frac{2q}{1+q^2}+2(1+q^2)(1+q^3)q^2
\left(-q^2+\frac{1}{1-q^6}+\frac{q^4}{1-q^8}+\frac{q^{16}}{(1-q^6)(1-q^8)}\right)\nonumber\\[3pt]
&\quad+2\sum_{k=3}^\infty (-q^2;q)_{2k-2}\,q^k \left(1-q^2+q^4+\sum_{n=0}^\infty b_nq^{2n}\right)\nonumber\\[3pt]
&=1-q-2q^4-2q^6-2q^7-2q^9+\frac{2q}{1+q^2}+\frac{2q^{18}+2q^{20}+2q^{21}+2q^{23}}
{(1-q^6)(1-q^8)}\nonumber\\[3pt]
&\quad +2(1+q^2)(1+q^3)q^2
\left(\frac{1}{1-q^6}+\frac{q^4}{1-q^8}\right)\nonumber\\[3pt]
&\quad +2\sum_{k=3}^\infty (-q^2;q)_{2k-2}q^k \sum_{n=0}^\infty b_nq^{2n}+2\sum_{k=3}^\infty (-q^3;q)_{2k-3}\,q^k (1+q^6).
\end{align}
Setting $a=0$, $b=10$ and replace $q$ with $q^2$ in Lemma \ref{lem-1}, we find that for $n\geq 33$, the coefficient of $q^n$ in
\[\frac{2q}{1+q^2}+\frac{2q^{21}}{(1-q^6)(1-q^8)}\]
is nonnegative. Thus the coefficient of $q^n$ in \eqref{equ-gf-n2mn-0-1} is nonnegative for $n\geq 33$, which implies that $\overline{N2}(0,n)\geq \overline{N2}(0,n-1)$
for $n\geq 33$. It is trivial to check that for $1\leq n\leq 32$, $\overline{N2}(0,n)\geq \overline{N2}(0,n-1)$ also holds. This yields \eqref{ine-nmn-ge-nmn-1-1-2} for $m=0$.

We proceed to show that \eqref{ine-nmn-ge-nmn-1-1-2} holds for $m\geq 1$. From Theorem \ref{lem-chan-mao-two-lem-1} and \eqref{equ-fg-n2-0-1}, we have
\begin{align}\label{equ-fg-n2-1-3-2-ge}
&\sum_{n=1}^\infty \left(\,\overline{N2}(m,n)-\overline{N2}(m,n-1)\,\right) q^n\nonumber\\[3pt]
=&2q f_{m,1}(q^2)+2(-q^2;q)_4 q^3f_{m,3}(q^2)+2\sum_{k=2\atop k\neq 3}^\infty(-q^2;q)_{2k-2}\,q^k f_{m,k}(q^2)\nonumber\\[3pt]
=&\frac{2q^{2m+1}}{1+q^2}+2(-q^2;q)_4 q^3f_{m,3}(q^2)+2\sum_{k=2\atop k\neq 3}^\infty(-q^2;q)_{2k-2}q^k f_{m,k}(q^2).
\end{align}
From Lemma \ref{lem-8}, we see that
 \begin{equation}
 f_{m,3}(q)=\sum_{n=-\infty}^\infty f_{n,2}(q)\,\frac{q^{3\mid m-n\mid}}{1-q^6}
 =f_{m,2}(q)+f_{m,2}(q)\,\frac{q^6}{1-q^6}+\sum_{n=-\infty\atop n\neq m}^\infty f_{n,2}(q)\,\frac{q^{3\mid m-n\mid}}{1-q^6}.
 \end{equation}
By Theorem \ref{lem-chan-mao-two-lem-2}, the coefficient of $q^n$ in $f_{m,2}(q)$ is nonnegative for all integer $m$ and $n\geq 0$. This allows us to transform $f_{m,3}(q)$ as follows:
 \begin{align}
 f_{m,3}(q)&=
 f_{m,2}(q)+\sum_{n=0}^\infty b_nq^n\nonumber\\[3pt]
 &=q^{m}\left(\frac{1-q^{ m+1}}{(1-q^2)(1-q^3)}+\frac{q^{ m+3}}{(1-q^3)(1-q^4)}\right)+\sum_{n=0}^\infty b_nq^n.
 \end{align}
 Hence
 \begin{align}\label{equ-exp-fm3-0}
 &2(-q^2;q)_4\,q^3f_{m,3}(q^2)\nonumber\\[3pt]
 =&2(-q^2;q)_4\,q^{2m+3}\left(\frac{1-q^{2m+2}}{(1-q^4)(1-q^6)}
 +\frac{q^{ 2m+6}}{(1-q^6)(1-q^8)}\right)+\sum_{n=0}^\infty b_nq^{n}
 \nonumber\\[3pt]
 =&2(-q^4;q)_2\,q^{2m+3}\frac{1-q^{2m+2}}{(1-q^2)(1-q^3)}
 +2(1+q^2)(1+q^5)\frac{q^{4m+9}}{(1-q^3)(1-q^4)}+
 \sum_{n=0}^\infty b_nq^{n}\nonumber\\[3pt]
 =&2q^{2m+3}\frac{1-q^{2m+2}}{(1-q^2)(1-q^3)}+2\frac{q^{4m+9}}{(1-q^3)(1-q^4)}+2(q^4+q^5+q^9)\left(q^{2m+3}\frac{1-q^{2m+2}}{(1-q^2)(1-q^3)}\right)\nonumber\\[3pt]
&+2(q^2+q^5+q^7)\left(\frac{q^{4m+9}}{(1-q^3)(1-q^4)}
\right)+\sum_{n=0}^\infty b_nq^{n}.
 \end{align}
 From Lemma \ref{lem-10}, we see that
 \[\frac{1-q^{2m+2}}{(1-q^2)(1-q^3)}\]
 has nonnegative coefficients.
 Together with \eqref{equ-exp-fm3-0}, we deduce that
 \begin{equation}\label{equ-exp-fm3}
 2(-q^2;q)_4\,q^3f_{m,3}(q^2)=2q^{2m+3}\frac{1-q^{2m+2}}{(1-q^2)(1-q^3)}
 +\frac{2q^{4m+9}}{(1-q^3)(1-q^4)}+
 \sum_{n=0}^\infty b_n q^n.
 \end{equation}
 Moreover, from Theorem \ref{lem-chan-mao-two-lem-2}, we see that
 \begin{equation}\label{equ-sum-k-ge-2-ne-3}
 \sum_{k=2\atop k\neq 3}^\infty(-q^2;q)_{2k-2}q^k f_{m,k}(q^2)=\sum_{n=0}^\infty b_nq^n.
 \end{equation}
Next we show that $\overline{N2}(m,n)\geq \overline{N2}(m,n-1)$ for $m\geq 2$.
Substituting \eqref{equ-exp-fm3} and \eqref{equ-sum-k-ge-2-ne-3} into \eqref{equ-fg-n2-1-3-2-ge}, we derive that
\begin{align}\label{equ-fg-n2-1-3-ex-f-1}
&\sum_{n=1}^\infty \left(\,\overline{N2}(m,n)-\overline{N2}(m,n-1)\,\right) q^n\nonumber\\[3pt]
=&\,\frac{2q^{2m+1}}{1+q^2}+2q^{2m+3}\frac{1-q^{2m+2}}{(1-q^2)(1-q^3)}
+\frac{2q^{4m+9}}{(1-q^3)(1-q^4)}
 +\sum_{n=0}^\infty b_n q^n\nonumber\\[3pt]
 =&\,\frac{2q^{2m+1}}{1+q^2}+2q^{2m+3}\frac{1-q^3+q^3-q^{2m+2}}{(1-q^2)(1-q^3)}
 +\frac{2q^{4m+9}}{(1-q^3)(1-q^4)}
 +\sum_{n=0}^\infty b_n q^n\nonumber\\[3pt]
 =&\,\frac{2q^{2m+1}}{1+q^2}+\frac{2q^{2m+3}}{1-q^2}
 +2q^{2m+6}\frac{1-q^{2m-1}}{(1-q^2)(1-q^3)}+\frac{2q^{4m+9}}{(1-q^3)(1-q^4)}
 +\sum_{n=0}^\infty b_n q^n\nonumber\\[3pt]
 =&\,\frac{2q^{2m+1}+2q^{2m+5}}{1-q^4}+2q^{2m+6}\frac{1-q^{2m-1}}{(1-q^2)(1-q^3)}
 +\frac{2q^{4m+9}}{(1-q^3)(1-q^4)}
 +\sum_{n=0}^\infty b_n q^n.
\end{align}
By Lemma \ref{lem-10}, we see that when $m\geq 2$,
\[q^{2m+6}\frac{1-q^{2m-1}}{(1-q^2)(1-q^3)}=\sum_{n=0}^\infty b_nq^n.\]
This gives $\overline{N2}(m,n)\geq \overline{N2}(m,n-1)$, as desired.

Finally, we consider the case $m=1$. In this case, by \eqref{equ-exp-fm3},
\begin{equation}\label{equ-exp-f13}
2(-q^2;q)_4\,q^3f_{1,3}(q^2)=2q^{5}\frac{1+q^{2}}{1-q^3}
 +\frac{2q^{13}}{(1-q^3)(1-q^4)}+
 \sum_{n=0}^\infty b_n q^n.
\end{equation}
Substituting \eqref{equ-sum-k-ge-2-ne-3} and \eqref{equ-exp-f13} into \eqref{equ-fg-n2-1-3-2-ge}, we see that
\begin{equation}
\sum_{n=1}^\infty \left(\,\overline{N2}(1,n)-\overline{N2}(1,n-1)\,\right) q^n
=\frac{2q^{3}}{1+q^2}+2q^{5}\frac{1+q^{2}}{1-q^3}
 +\frac{2q^{13}}{(1-q^3)(1-q^4)}+
 \sum_{n=0}^\infty b_n q^n.
\end{equation}
From Lemma \ref{lem-1}, we find that for $n\geq 19$, the coefficient of $q^n$ in
\[\frac{2q^{3}}{1+q^2}+\frac{2q^{13}}{(1-q^3)(1-q^4)}\]
is nonnegative. This gives $\overline{N2}(1,n)\geq \overline{N2}(1,n-1)$ for $n\geq 19$. It can be checked that for $1\leq n\leq 18$,  $\overline{N2}(1,n)\geq \overline{N2}(1,n-1)$ still holds. This completes the entire proof.\qed

\section{The proof of Theorem \ref{main-thm-2}} \label{sec:4}

In this section, we give a proof of Theorem \ref{main-thm-2}. To this end, we need the following lemma.

\begin{lem}\label{lem-2}
For integer $k\geq 0$, let
\[\frac{1}{(qz;q)_k\,(q/z;q)_k}=\sum_{n=0}^\infty \sum_{m=-\infty}^\infty a_{k,m}(n)z^mq^n.\]
Then for $m\geq 0$, we have 
$a_{k,m}(n)\geq a_{k,m+2}(n)$.
Equivalently, for $m\geq 0$, the coefficient of $z^m q^n$ in
\[\frac{1-z^{-2}}{(qz;q)_k(q/z;q)_k}\]
is nonnegative.
\end{lem}

\pf By definition, we see that
\begin{equation}\label{sym-akmn}
a_{k,m}(n)=a_{k,-m}(n).
\end{equation}
Moreover, it is clear that
\begin{align}
\sum_{n=0}^\infty \sum_{m=-\infty}^\infty a_{k+1,m}(n)z^mq^n&=\frac{1}{(qz;q)_{k+1}(q/z;q)_{k+1}}\nonumber\\[3pt]
&=\frac{1}{(1-zq^{k+1})(1-q^{k+1}/z)}\sum_{n=0}^\infty \sum_{m=-\infty}^\infty a_{k,m}(n)z^mq^n\nonumber\\[3pt]
&=\sum_{r=0}^\infty \sum_{i=0}^r z^{r-2i}q^{r(k+1)}\sum_{n=0}^\infty \sum_{m=-\infty}^\infty a_{k,m}(n)z^mq^n.
\end{align}
Thus we have
\begin{equation}\label{equ-ref-akmn}
a_{k+1,m}(n)=\sum_{r=0}^{\lfloor \frac{n}{k+1} \rfloor}\sum_{i=0}^r a_{k,m-r+2i}\left(n-r(k+1)\right).
\end{equation}
 
We prove this lemma by induction on $k$. For $k=1$, it is trivial to check that
\[a_{1,m}(n)=\begin{cases}
1& \text{if } m\equiv n\pmod{2}\text{ and } n\geq |m|;\\
0& \text{otherwise.}
\end{cases}\]
This gives our desired result.

Set $b_{k,m}(n)=a_{k,m}(n)-a_{k,m+2}(n)$ and assume that  $b_{k,m}(n)\geq 0$ for $m\geq 0$. From \eqref{equ-ref-akmn}, we derive that
\begin{equation}\label{equ-ref-bkmn}
b_{k+1,m}(n)=\sum_{r=0}^{\lfloor \frac{n}{k+1} \rfloor} \sum_{i=0}^r b_{k,m-r+2i}\left(n-r(k+1)\right).
\end{equation}
Moreover, by \eqref{sym-akmn}, we see that
\begin{equation}\label{equ-sym-bkmn-1}
b_{k,m}(n)=-b_{k,-m-2}(n)
\end{equation}
and therefore 
\begin{equation}\label{equ-sym-bkmn}
\sum_{r=m+1}^{\lfloor \frac{n}{k+1} \rfloor} \sum_{i=0}^{r-m-1} b_{k,m-r+2i}\left(n-r(k+1)\right)=0.
\end{equation}
Thus by \eqref{equ-ref-bkmn} and \eqref{equ-sym-bkmn}, we derive that for $m\geq 0$,
\begin{align}\label{equ-lem-bk+1mn}
b_{k+1,m}(n)&=\sum_{r=0}^m \sum_{i=0}^r b_{k,m-r+2i}(n-r(k+1))+\sum_{r=m+1}^{\lfloor \frac{n}{k+1} \rfloor} \sum_{i=0}^r b_{k,m-r+2i}(n-r(k+1))\nonumber\\[3pt]
&=\sum_{r=0}^m \sum_{i=0}^r b_{k,m-r+2i}(n-r(k+1))+\sum_{r=m+1}^{\lfloor \frac{n}{k+1} \rfloor} \sum_{i=r-m}^r b_{k,m-r+2i}(n-r(k+1)).
\end{align}
From induction hypothesis, we find that each term in the above summation is nonnegative. Thus $b_{k+1,m}(n)\geq 0$.
This completes the proof.\qed

We now give a proof of Theorem \ref{main-thm-2}.

\noindent{\it Proof of Theorem \ref{main-thm-2}.} By \eqref{equ-gf-over-rank}, for $m\geq 0$, $\overline{N}(m,n)\geq \overline{N}(m+2,n)$ is equivalent to that the coefficient of $z^m$ in
\begin{equation}\label{equ-gf-over-rank-1}
\sum_{k=0}^\infty\frac{(-1;q)_k \, q^{k(k+1)/2}(1-z^{-2})}{(zq;q)_k(q/z;q)_k}
\end{equation}
is nonnegative.
But by Lemma \ref{lem-2},
\begin{equation}
[z^m]\,\sum_{k=0}^\infty\frac{(-1;q)_k \, q^{k(k+1)/2}(1-z^{-2})}{(zq;q)_k(q/z;q)_k}
=\sum_{k=0}^\infty(-1;q)_k \, q^{k(k+1)/2}[z^m]\,\frac{1-z^{-2}}{(zq;q)_k(q/z;q)_k},
\end{equation}
which is clearly has nonnegative coefficients, where $[z^m]\,f(z)$ denotes the coefficient of $z^m$ in $f(z)$. This yields \eqref{ine-nmn-ge-nmn-1-1-2-1}.

Similarly, by \eqref{equ-gf-over-m2-rank}, for $m\geq 0$, $\overline{N2}(m,n)\geq \overline{N2}(m+2,n)$ is equivalent to that the coefficient of $z^m$ in
\begin{equation*}
\sum_{k=0}^\infty\frac{(-1;q)_{2k}\,q^{k}(1-z^{-2})}{(zq^2;q^2)_k\,(q^2/z;q^2)_k}
\end{equation*}
is nonnegative. Again using Lemma \ref{lem-2}, we see that
\begin{equation*}
[z^m]\,\sum_{k=0}^\infty\frac{(-1;q)_{2k}\,q^{k}(1-z^{-2})}{(zq^2;q^2)_k\,(q^2/z;q^2)_k}
=\sum_{k=0}^\infty(-1;q)_{2k}\,q^{k}[z^m]\,\frac{(1-z^{-2})}{(zq^2;q^2)_k\,(q^2/z;q^2)_k},
\end{equation*}
which has nonnegative coefficients. This completes the proof.\qed

\section{Conclusions}
The rank of partitions gives combinatorial interpretations of several Ramanujan's famous congruence formulas.
In this paper, we derive several monotonicity inequalities of the $D$-rank and $M_2$-rank for overpartitions and use them to prove a conjecture of Chan and Mao \cite{Chan-Mao-2014}. 
Our proofs are based on the study of generating functions for such ranks of overpartitions, which are analytic.  It would be interesting to find bijective proofs for our results. We will work on this in the future.

\section*{Acknowledgments}
 The first author acknowledges support from the Swiss National Science Foundation (Grant number P2ZHP2\_171879).
This work was done during the first author's visit to the Harbin Institute of Technology (HIT). The first author would like to thank  Prof. Quanhua Xu and the second author for the hospitality.


\end{document}